\documentclass[10pt]{article}
\usepackage[cp866]{inputenc}
\usepackage{amsmath,amsfonts,amssymb}
\begin{document}

\hfill MSC 20B20,20B22

\begin{center}
{\Large On correspondence between right near-domains
and sharply 2--transitive groups}\\
{\bf Andrey A. Simonov}\\
a.simonov@g.nsu.ru
\end{center}

\begin{quote}
\noindent{\sc Abstract. } The right near-domain
is defined to loosen near-domain axioms.
Correspondence of a class of the right near-domains and a class of
sharply 2--transitive groups is constructed.
\medskip

\end{quote}

In \cite {kar1, kar2} for exposition of
{\em sharply 2--transitive groups} the concept {\em near-domain}
is introduced. Near-domain is an
algebraic system $ (B, \cdot, +, {}^{-1}, 0) $
with two binary operations $ \cdot, + $
for which axioms hold:

1. $(B, +)$ is a loop with a unit element $0$;

2. $a + b = 0\Rightarrow b+a=0$;

3. $(B_1, \cdot, {}^{-1})$ is a group with an unit element $e$, where $B_1 = B
\setminus \{0\}$;

4. $(\forall x \in B)\quad  x\cdot 0= 0$;

5. $(\forall x,y,z \in B) \quad  (x+y)\cdot z = x\cdot z+ y\cdot z$;

6. $(\forall a,b \in B) (\exists $ $r_{a,b}\in B_1) \quad
(x+a) + b = x\cdot r_{a,b}+ (a + b)$ for any $x\in B$.
\\
Until recently it is not known any example of a near-domain which is not
a near-field. In the given work it is offered to loosen near-domain axioms, having
left only necessary ones for construction
of sharply 2--transitive groups. In particular,
it is offered to refuse from axioms 2, 4 and to loosen axioms 1, 5.

Let's define the right near-domain as algebraic system $ (B, \cdot, +, -,
 {}^{-1}, 0) $ with operations:
$
(+):B\times B_1\to B,$  $(-):B\times B_1\to B,$ $(\cdot ):B\times B_1\to B$,
 where  $B_1 = B \setminus \{0\}$,
for which axioms are hold
\\
A1. $(\forall x\in B)(\forall y\in B_1) \ (x-y)+y=x$;
\\
A2. $(\forall x\in B)(\forall y\in B_1) \ (x+y)-y=x$;
\\
A3. $(\forall x\in B_1) \ $ $x-x=0$;
\\
A4. $(B_1, \cdot , {}^{-1})$ is a group with an unit element $e\in B_1$;
\\
A5. $(\forall x\in B)(\forall y, z\in B_1)(\exists  \ h(y,z)\in B_1)
\ (x+y)z=xh(y,z)+yz$;
\\
A6. $(\forall x\in B)(\forall y, z\in B_1: \
y+z\neq 0)(\exists  \ r(y,z)\in B_1)$
$(x+y)+z=xr(y,z)+(y+z)$;
\\
A7. $(\forall x\in B)(\forall z\in B_1)
(\exists  \ v(z)\in B_1) \ (x+(0-z))+z=xv(z)$.

Axioms A1---A3 define algebraic system $ (B, +,-, 0) $ as the right loop.
We define labels $L (x) =0-x $ then from A1 follows $L (x) +x=0$.
Thus the map $L:B_1\to B_1$ defines left inverse in the right loop.

Let's consider now the elementary consequences of axioms.\\
{\bf Lemma.} {\em In the right near-domain the following properties hold:\\
1. $(\forall x\in B_1) \ 0x=0$;\\
2. $h(x,y)=EL(x)L(xy)$, where $E(x)=x^{-1}$,
$EL$ --- superposition of transformations $L$ and $E$;\\
3. $r(y,z)=E(L(z)-y)L(y+z)$;\\
4. $x-z=xEv(z)+L(z)$;\\
5. $v(z)= EL^{2}(z)z$
}

Let's define a map $u:B_1\to B $ by the rule $u(x)=0x$.

From A5 follows, that
$(\forall x, y \in B_1) \
(L(x)+x)y=L(x)h(x,y)+xy=u(y)$, hence
\begin{equation}
\label{h1}
h(x,y)=EL(x)(u(y)-xy).
\end{equation}
If we will sequentially apply A5 for arbitrary $z, t \in B_1$  then
we receive:
$$
h(y,z)h(yz,t)=h(y,zt).
$$
Let's write the given equality applying the
expression (\ref {h1}). With the reduction account,
we will receive equality $: (u (z)-yz) EL (yz) =e $, hence,
$u (z) =L (yz) +yz=0$.
Thus, the first and second conditions of the lemma are satisfied.

Let's consider now consequences from A6. Let
$x=L(y+z)(r(y,z))^{-1}$ $\Rightarrow $
$(L(y+z)(r(y,z))^{-1}+y)+z=0$, whence we will receive the expression from
the third condition of the lemma.

In the case when $y+z=0$, we will consider consequences from
A7 and A2: $x+L (z) =xv (z)-z $.
We define $x'=xv(z) $,
hence the fourth condition of the lemma is fulfilled
$x'E(v (z))+L (z) =x'-z $.

Let's note A2  with the condition of the received expression
$(x+z)-z = (x+z) Ev(z) +L (z) =x $.
At $x=0$ we will receive equality $zEv(z) +L (z) =0$. Then with the account
$L^2 (x) =LL (x) $, we will come to justice of
the fifth condition of the lemma.
\hfill $\Box $

The operation "$ - $" is expressed through the
operations "$ + $", "$ \cdot $", $L $, $E $ hence
we will understand algebraic system $ (B, \cdot, +, - ,{}^{-1} , 0) $
as $(B, \cdot , +, {}^{-1}, L, 0)$.
\\

Let's consider the algebraic system
$ (H, \cdot , \phi ,{}^{-1}, 0) $ from \cite{sim1},
with the operations:
$$
(\cdot): H\times H_1\to H, \phi :H\to H, \mbox{ where } H_1=H\setminus \{0\},
$$
for which the following axioms are fulfilled:

F1. $(H_1, \cdot ,{}^{-1})$ is a group with an unit element $e$;

F2. $0 x = 0, x\in H_1$;

F3. $\phi (e) = 0$;

F4. $\phi (\phi (x)\phi (y))=
\phi (x\phi (y^{-1}))y, x\in H, y\in H_1\setminus \{e_1\}$,
\\
The similar algebraic system was investigated in \cite{les}.
\\
{\bf Theorem 1.}
{\em The class of algebraic systems $ (B, \cdot, +,{}^{-1}, L, 0) $ and
$ (B, \cdot,{}^{-1}, \phi, 0) $ are rational equivalent.}
\\
Let's introduce a map $ \phi:B\rightarrow B $, defined in an aspect
$ \phi (x) =x (0-e) +e=xa+e $.
Let's calculate quadrate of function $ \phi $ taking into account
the conditions two and five of the lemma:
$$
\phi ^2(x)=(xa+e)a+e=(xL(a)+a)+e=xL(a)EL^2(e)=x.
$$
From the definition follows
$\phi (e)=a+e=0$ and $\phi (0)=e$.
By means of the map $ \phi $ it is possible to express additive operation.
Really, $ \phi (x) y = (xa+e) y=xL (y) +y $, hence,
if $x=zEL (y) $, then $z+y =\varphi (zEL (y)) y $.
Let's rewrite now identity from A2:
$z=(z+y)-y=\phi (zEL(y))y-y$. Having introduced labels $t =\phi (zEL (y)) y $, we
express $z=\phi (ty^{-1})L(y)$, then \ $t-y=\phi (ty^{-1})L(y)$.

Calculating the
value $t = (x+z) - (y+z) $ in the case $y\neq L (z) $, using at first
A2: $ (x+z) =t + (y+z) $, and then the third identity of
the lemma:
$(x+z)=(t(r(y,z))^{-1}+y)+z$. Applying twice identity from
A2, we have the identity:
$$
(x+z)-(y+z)=(x-y)(L(z)-y)^{-1}L(y+z).
$$

Let's rewrite the given identity with the account $y\neq e, z=L ^ {-1} (e) $
replacing additive binary operations by their expressions through
the function $\phi $:
\begin{equation}
\label{phi1}
\phi (\phi (x)E\phi (y)) =\phi (xy^{-1})E\phi E
(y)=\phi (xy^{-1})E\phi E(y).
\end{equation}
At $x=0$ this identity takes a simple form $ \phi E\phi (y) =E\phi E (y) $,
using it, we note identity (\ref {phi1}) for  $y=E\phi E(t)$:
\begin{equation}
\label{varphi}
\phi (\phi (x)\phi (t)) =\phi (x\phi E(t))t.
\end{equation}
Thus, we have the map
$\mathbb A:(B, \cdot, +,{}^{-1}, L, 0) \to (B, \cdot,{}^{-1}, \phi, 0)$.

Let's make the inverse construction. We will consider expression from F4 at
$x=e, y=t ^ {-1} $,
then under condition F2 and F3 we come to equality
$ \varphi ^2 (t) = \varphi (0) t $.
On one hand $ \varphi ^4 (t) = (\varphi (0)) ^ {2} t $, and on the other hand
$ \varphi ^4 (t) = \varphi (\varphi ^2 (\varphi (t))) = \varphi (\varphi (0) \varphi (t)) $.
It is also possible to note the
last expression with the account F4 and F2:
$ \varphi (\varphi (0) \varphi (t)) = \varphi (0\varphi (t ^ {-1})) t =\varphi (0) t $.
Thus, we come to equality $ \varphi ^2 (0) = \varphi (0) $, hence,
$ \varphi (0) =e $ and $ \varphi ^2 (t) =t $.

From F4 for $x=E \varphi E(y)$ follows, that $(\forall y\in B_1\setminus \{e\})
\ \varphi E \varphi (y)=E\varphi E(y) $.

By means of arbitrary bijection $L: B_1\to B_1$ we introduce operations
$$
x+y=\varphi (xEL(y))y, \ x-y=\varphi (xy^{-1})L(y).
$$
With the account of F2, F3 and $ \varphi ^2=id $ it is easy to check up
the performance of the axioms A1---A3
of the right loop. The performance A5 follows from the operation definition:
$$
(x+y)z=\varphi (xEL(y))yz=\varphi (xEL(y)L(yz)EL(yz))yz=xEL(y)L(yz)+yz.
$$
Then we take advantage of identities $\varphi ^2=id$,
$\varphi E \varphi =E\varphi E$ and F4 to receive A6:	
$$
(x+y)+z=\varphi (\varphi (xEL(y))yEL(z))z =
$$
$$
\varphi (xEL(y)\varphi E\varphi (yEL(z)))\varphi (yEL(z))z=
$$
$$
\varphi (xEL(y)\varphi E\varphi (yEL(z))
L[\varphi (yEL(z))z]
EL[\varphi (yEL(z))z])
\varphi (yEL(z))z=
$$
$$
xEL(y)\varphi E\varphi (yEL(z))L[\varphi (yEL(z))z]+(y+z)=
$$
$$
xE(\varphi (L(z)E(y))L(y))L[\varphi (yEL(z))z]+(y+z)=
xE(L(z)-y)L(y+z)+(y+z).
$$
Now we take advantage of identity $\varphi ^2=id$  for construction of
expression A7:
$$
(x+L(z))+z=\varphi (\varphi (xEL^2(z))L(z)EL(z))z=
\varphi ^2(xEL^2(z))z=xEL^2(z)z.
$$
For any bijection $L $ we have constructed the map
$\mathbb F_L: (B, \cdot,{}^{-1}, \varphi, 0) \to (B, \cdot, +,{}^{-1},$
$L, 0)$
so, that the algebraic systems $(B, \cdot,{}^{-1}, \varphi, 0)$ and
$\mathbb A \circ \mathbb F_L (B, \cdot,{}^{-1}, \varphi, 0)$ are isomorphic.
In the opposite direction the algebraic systems
$(B, \cdot, +',{}^{-1}, L', 0)$ and
$ \mathbb F_L \circ \mathbb A (B_1, 0, L ', \cdot, + ')$ are isomorphic
only at $L=L ' $.
\hfill $\Box $

The group $T_2 (B) $ of transformations of a set $B$
is called sharply 2--transitive group,
if for arbitrary pairs $ (x_1, x_2) \neq (y_1, y_2) \in \widehat {B^2} $, where
$ \widehat {B^2} =B^2 \setminus \{(x, x) |x\in B \} $ there exists
an unique element
$g\in T_2 (B) $ for which the equalities $g (x_1) =y_1$ and $g (x_2) =y_2$ are 
held.

{\bf Theorem 2.} {\em The class
of algebraic systems $ (B, \cdot,{}^{-1}, \varphi, 0) $
and the class of sharply 2--transitive groups $T_2 (B) $ are
rational equivalent.} \\ On the set
$ \widehat {B^2} $ we define a function $f:B\times \widehat {B^2} \to B $
as
\begin{equation}
\label{f}
f(x,y_1,y_2)=\varphi (x \varphi (y_1 y_2^{-1}))y_2,
\end{equation}
if $y_2 \neq 0$ and $f (x, y_1,0) =x y_1$ otherwise.
Not to consider two cases separately, we, by means of
multiplicative partial operation $(\cdot ): B\times B_1\to B$,
define the groupoid on $B $ so, that  $(\forall x \in B) \
x 0 = \varphi (x)$, $0^{-1}=0$.

Let's define a binary operation $G $ on
the set $ \widehat {B^2} $ in the form
\begin{equation}
\label{G}
\left(
\begin{array}{c}
x_1 \\
x_2
\end{array}
\right) \left(
\begin{array}{c}
y_1 \\
y_2
\end{array}
\right) =
\left(
\begin{array}{c}
f(x_1,y_1,y_2) \\
f(x_2,y_1,y_2)
\end{array}
\right) =
\left(
\begin{array}{c}
\varphi (x_1\varphi (y_1y_2^{-1}))y_2 \\
\varphi (x_2\varphi (y_1y_2^{-1}))y_2
\end{array}
\right) .
\end{equation}
Supposing, that there are
pairs $ (x_1, x_2) \neq (y_1, y_2) \in \widehat {B^2} $,
that $f (x_1, y_1, y_2) =f (x_2, y_1, y_2) $. Then, for $y_2\neq 0$ after
multiplication of
the both parts of equality on $ y_2 ^ {-1} $ and transformations
by the function $ \varphi $, we will come to equality
$x_1\varphi (y_1y_2 ^ {-1}) =x_2\varphi (y_1y_2 ^ {-1}) $ from which follows,
that $x_1=x_2$. At $y_2=0$ we get
the equality $x_1y_1=x_2y_1$, hence, $x_1=x_2$.
We have come to an inconsistency.
Thus, the operation $G $, defined above,  is a magma.

It is easy to check, that the pair $ (e, 0) \in \widehat {B^2} $
is the left and the right unit element.
Now we check the associativity:
$$
\varphi (\varphi (x_i\varphi (y_1y_2^{-1}))y_2\varphi
(z_1z_2^{-1})) z_2=\varphi ( \varphi (x_i\varphi
(y_1y_2^{-1}))\varphi \varphi ( y_2\varphi (z_1z_2^{-1})))
z_2=
$$
$$\varphi (x_i\varphi (y_1y_2^{-1})\varphi E\varphi (y_2\varphi
(z_1z_2^{-1}))) \varphi (y_2\varphi (z_1z_2^{-1}))
z_2=
$$
$$
\varphi (x_i\varphi (y_1\varphi (z_1z_2^{-1})E\varphi
(y_2\varphi (z_1z_2^{-1}))) \varphi (y_2\varphi
(z_1z_2^{-1})) z_2.
$$
We
have come to a semigroup with a unit element. We will discover now the left inverse:
$$
\left(
\begin{array}{c}
x_1 \\
x_2
\end{array}
\right) ^{-1}\left(
\begin{array}{c}
x_1 \\
x_2
\end{array}
\right) =
\left(
\begin{array}{c}
\varphi (x_2^{-1})E\varphi (x_1x_2^{-1}) \\
E\varphi (x_1x_2^{-1})
\end{array}
\right)
\left(
\begin{array}{c}
x_1 \\
x_2
\end{array}
\right) =
\left(
\begin{array}{c}
e \\
0
\end{array}
\right).
$$
Let's check that it is also the right inverse:
$$
\varphi (x_i\varphi (\varphi (x_2^{-1})E\varphi (x_1x_2^{-1})\varphi
(x_1x_2^{-1})) )E\varphi (x_1x_2^{-1})=
\varphi (x_ix_2^{-1})E\varphi (x_1x_2^{-1}).
$$
Thus, we defined that $G$ is a group, but since it operates on the set
$ \widehat {B^2} $ sharply transitive, then
the group $G $ at an operation on the set $B $
it will be sharply 2--transitive, hence, we have constructed the map
${\mathbb G}:(B, \cdot,{}^{-1}, \varphi, 0)\to T_2(B)$.

Now we make the inverse construction and on a group $T_2 (B) $
we will construct an algebraic system $ (B, \cdot,{}^{-1}, \phi, 0) $.
For an arbitrary pair $ (e_1, e_2) \in \widehat {B^2} $ it is
possible to construct the bijective map $T_2 (B) \to \widehat {B^2} $,
putting in correspondence to an element $g\in T_2 (B) $ the pair
$ [x_1, x_2]$ so that
\begin{equation}
\label{e_i}
(e_1, e_2) \cdot g=
(e_1, e_2) \cdot [x_1, x_2]=
(e_1\cdot [x_1, x_2], e_2 \cdot [x_1, x_2])=
(x_1,x_2).
\end{equation}
The given bijection induces 
the isomorphic group $G \simeq  T_2 (B)$ on the set of pairs
$\widehat {B^2}$.
The pair $[e_1, e_2]$ is an unit of the group $G$.

At serial transformation of the pair $ (e_1, e_2)$ by elements $ [x_1, x_2] $
and $ [y_1, y_2] $ we come to equality:
\begin{equation}
\label{g}
[x_1,x_2][y_1,y_2]=[x_1\cdot [y_1,y_2],x_2\cdot [y_1,y_2]],
\end{equation}
from which, with the account (\ref {e_i}), follows, that on a subset
$B_1= \{x\in B|[x, e_2]\in G\}$
it is possible to introduce the group structure naturally.
The map $e_1\cdot [x, e_2] \mapsto x$ induces on $B_1$ a group structure.
Multiplication in the group $B_1$, as well as in the group $T_2 (B) $
we will write without a point. We will expand the group operation to a
partial operation $B\times B_1 \to B $, having predetermined it in an aspect
$e_2y=e_2\cdot [y, e_2] =e_2$ so, that $e_2$ will be the left zero in a
partial operation  $(\cdot ): B\times B_1 \to B$.

From (\ref {e_i}) and (\ref {g}) follows,$ [e_2, e_1] $
is an involution of $G$. We define
$ \phi:B\to B $ in an aspect $ \phi (x) = x\cdot [e_2, e_1] $, then
$\phi (e_1)=e_2$ and
\begin{equation}
\label{e_2,e_1}
[e_2, e_1][x_2, x_1]=[x_1, x_2]=[\phi (x_1), \phi (x_2)] [e_2, e_1].
\end{equation}
For an arbitrary
$ [e_1, x_2] \in G $, at $x_2 \in B_1 \setminus \{e_1 \}$ it is possible to note:
$$
[e_1,x_2]=[x_2^{-1},e_1][x_2,e_2]=[\phi (x_2^{-1}),e_2][e_2,e_1][x_2,e_2].
$$
On the other hand, with the account (\ref {e_2,e_1}) for $ [e_1, x_2] $ it is 
fair
$$
[e_1,x_2]=[e_2,e_1][\varphi (x_2),e_2][e_2,e_1].
$$
Having taken advantage of the two received expressions and equating outcomes
of transformation arbitrary $t\in B $ by element $ [e_1, x_2] \in G$, we
come to identity:
$$
\phi (\phi (t)\phi (x_2))=
\phi (t\phi (x_2^{-1}))x_2, \ t\in B, x_2\in B_1 \setminus \{e_1 \}.
$$
The map ${\mathbb F}_{(e_1,e_2)}:T_2(B)\to (B, \cdot,{}^{-1}, \phi, 0)$
is constructed,
putting in correspondence to group $T_2 (B) $ algebraic system
$(B, \cdot,{}^{-1}, \phi, 0)$.

Let's notice still, that for arbitrary $ [x_1, x_2] \in T_2 (B) $ it is 
possible to note:
$$
[x_1,x_2]=
\left\{
\begin{array}{ll}
[\phi (x_1x_2^{-1}), e_2] [e_2, e_1] [x_2, e_2] ,
&  x_2\in B_1,\\
\left[
x_1
, e_2\right]
, &
x_2=e_2.
\end{array}
\right.
$$
Then for arbitrary $t\in B$ under condition of $x_2\neq e_2$ and $t\cdot [x_1,e_2]=tx_1$
the equality:
\begin{equation}
\label{ff}
t\cdot [x_1,x_2]=t\cdot [\phi (x_1x_2^{-1}), e_2] [e_2, e_1] [x_2, e_2]=
\phi (t\phi (x_1x_2^{-1})) x_2
\end{equation}
is fair.
Comparing (\ref {f}), (\ref {G}) with (\ref {ff}) and (\ref {g}) we
come to that there is a natural isomorphism ${\mathbb G} \circ {\mathbb F}_{(e_1,e_2)}: T_2(B)
\to T_2'(B)$, thus ${\mathbb G} \circ {\mathbb F}_{(e_1,e_2)} = id$.
Isomorphism of algebraic systems
${\mathbb F}_{(e_1,e_2)} \circ {\mathbb G}:
(B, \cdot,{}^{-1}, \phi, 0)\to (B', \cdot ',{}^{-1}, \phi, e_2)$
is set by map
${\mathbb F}_{(e_1,e_2)} \circ {\mathbb G}:
x\mapsto \varphi (x\varphi (e_1e_2^{-1})) e_2$, thus
${\mathbb F}_{(e_1,e_2)} \circ {\mathbb G} = id$.
\hfill $\Box $

Let's consider some examples of the right near-domains constructed over
a skew field $ \mathbb K $ for which $ \varphi (x) =-x+1$,
$x\in \mathbb K $. As the first example we consider $L (x) =ax $:
$$
x\oplus y=-xa^{-1}+y, \
x\ominus y=-xa+ay, \ r(y,z)=-a^{-1}, \ v(z)=a^{-2}.
$$
In such right near-domain bilaterial distributivity is fulfilled and the
identity $L (x\oplus y) = L (x) \oplus L (y) $ is hold.
For the second example over a skew field we consider $L (x) =-x^{-1} $, then
$$
x\oplus y=xy^2+y,\
x\ominus y=xy^{-2}-y^{-1},\
r(y,z)=y^2z(z+y)^{-1}(yz+1), \ h(y,z)=z^{-1}.
$$
For the given right loop $L(x\oplus y)\neq L(x)\oplus L(y)$, but it is fulfilled
$L(x)\oplus x=x\oplus L(x)=0$.


\begin{thebibliography}{9}

\bibitem{kar1}
{\em  Karzel H.} Inzidenzgruppen I. Lecture Notes by Pieper, I. and Sorensen, K., University of
Hamburg (1965), 123-135.

\bibitem{kar2}
{\em  Karzel H.} Zusammenhange zwischen Fastbereichen, scharf zweifach transitiven
Permutationsgruppen und 2-Strukturen mit Rechtecksaxiom, Abh. Math. Sem. Univ. Hamburg
32 (1968), 191-206.

\bibitem{sim1}
{\em  Simonov A.A.}
About correspondence between neardomains and groups. Algebra i Logic.
vol. 45, 2,
2006.

\bibitem{les}
{\em  Leissner W.}
On the Functional Equation $\phi(xy^{-1}) =\phi(\phi(x)\phi(y)^{-1})\phi (y^{-1})
$ over a Group. Report of Meetings. Elfte internationale Tagung \"uber
Funktionalgleichungen in Oberwolfach vom 14. bis 20. Dezember 1973.

\bibitem{mal}
{\em  Maltsev A. I.}
Structural performance of some classes of algebras,
Doklady of the Academy of Sciences of the USSR, 120, No. 1, 29-32, 1958.


\end{thebibliography}
\end{document}